\input amstex
\documentstyle{amsppt}

\magnification=\magstep1
\hsize=6.5truein
\vsize=9truein

 at 10truept
 at 7truept
 at 10truept 
 at 10truept

\document

\baselineskip=.15truein

\def\g{{\frak{g}}}

\def\RM{{\Bbb{R}}}

\def \ot {\otimes}
\def \gh{\oplus\underline{\g^{\otimes \bullet}} }
\def \ghu{\oplus\underline{\g_1^{\otimes \bullet}} }
\def \ghd{\oplus\underline{\g_2^{\otimes \bullet}} }
\def \gu{\ve^\bullet\!\underline{\g_1^{\otimes \bullet}}}
\def \gq{\ve^\bullet\!\underline{\g^{\otimes \bullet}}}
\def \gd{\ve^\bullet\!\underline{\g_2^{\otimes \bullet}}}

\def \End {\mathop{\hbox{\rm End}}\nolimits}

\def \Hom {\mathop{\hbox{\rm Hom}}\nolimits}

\def \HKR {{\scriptstyle{\hbox{HKR}}}}

\def \Id {\mathop{\hbox{\rm id}}\nolimits}

\def \ve{\wedge}

\def \p {\partial}

\topmatter

\title  
Lifts of $C_\infty$ and $L_\infty$-morphisms to $G_\infty$-morphisms
\endtitle

\author
Gr\'egory Ginot${}^{(a)}$ \, , \;  Gilles Halbout${}^{(b)}$ 
\endauthor 

\leftheadtext{ Gregory Ginot, Gilles Halbout } 
\rightheadtext{ Lift of $C_\infty$ and $L_\infty$ morphisms to $G_\infty$ morphisms } 

\affil
${}^{(a)}$ \! Laboratoire Analyse G\'eom\'etrie et Applications, Universit\'e Paris~13 et Ecole Normale Sup\`erieure de Cachan,  France\\   
 ${}^{(b)}$ \! Institut de Recherche Math\'ematique Avanc\'ee,   
Universit\'e Louis Pasteur et CNRS, Strasbourg, France \\ 
\endaffil

\address\hskip-\parindent
${}^{(a)}$ \! Laboratoire Analyse G\'eom\'etrie et Applications, Universit\'e Paris~13, \newline   
Centre des Math\'ematiques et de Leurs Applications, ENS Cachan,\newline
61, av. du Pr\'esident Wilson, 94230 Cachan,  France \newline
e-mail:ginot\@{}cmla.ens-cachan.fr \newline
  {}  \newline
  ${}^{(b)}$ \! \hbox{Institut de Recherche Math\'ematique Avanc\'ee,}  \newline   
Universit\'e Louis Pasteur et CNRS\newline
7, rue Ren\'e Descartes, 67084 Strasbourg, France \newline 
e-mail:halbout\@{}math.u-strasbg.fr \newline
\endaddress

\abstract 
Let  $\g_2$ be the
Hochschild complex of cochains on $C^\infty(\RM^n)$ and $\g_1$ be the space of multivector fields on $\RM^n$. 
In this paper we prove that given any  $G_\infty$-structure ({\rm i.e.} 
Gerstenhaber algebra up to homotopy structure) on $\g_2$,
and any $C_\infty$-morphism $\varphi$ ({\rm i.e.} morphism of commutative, 
associative algebra up to homotopy) between
$\g_1$ and $\g_2$, there exists a $G_\infty$-morphism $\Phi$
between
$\g_1$ and $\g_2$ that restricts to $\varphi$. We also show that 
any $L_\infty$-morphism ({\rm i.e.} morphism of Lie algebra up to homotopy), 
in particular the one constructed by Kontsevich,  can be deformed into a $G_\infty$-morphism,  using Tamarkin's 
method for any $G_\infty$-structure on $\g_2$. We also show that 
any two of such $G_\infty$-morphisms are homotopic.
\endabstract

\endtopmatter

\footnote""{{\bf Keywords~:} Deformation quantization,
star-product, homotopy formulas, homological methods.}  
\footnote""{ {\bf AMS Classification~:} 
Primary 16E40, 53D55, Secondary 18D50, 16S80. }

\vskip20pt

\centerline {\bf 0-Introduction }

\vskip13pt

Let $M$ be a differential manifold and $\g_2=(C^\bullet(A,A),b)$ be the Hochschild cochain complex on $A=C^\infty(M)$. The classical Hochschild-Kostant-Rosenberg theorem states that the cohomology of $\g_2$ is the graded Lie algebra 
$\g_1=\Gamma(M,\ve^\bullet TM)$ of multivector fields on
$M$. There is also a graded Lie algebra structure on $\g_2$ given by the Gerstenhaber bracket. In particular $\g_1$ and $\g_2$ are also Lie algebras up to homotopy ($L_\infty$-algebra for short). In the case $M=\RM^n$,
using different methods, Kontsevich ([Ko1]
and [Ko2]) and Tamarkin ([Ta]) have proved
the existence of Lie homomorphisms ``up to homotopy'' ($L_\infty$-morphisms) from $\g_1$ to $\g_2$.
Kontsevich's proof uses graph complex and is related to multizeta functions whereas
Tamarkin's construction uses the existence of Drinfeld's associators.
In fact  Tamarkin's $L_\infty$-morphism
comes from the restriction of a Gerstenhaber algebra up to homotopy 
homomorphism ($G_\infty$-morphism) from $\g_1$ to $\g_2$. 
The $G_\infty$-algebra structure on $\g_1$ is induced by its classical 
Gerstenhaber algebra structure and  a far less trivial 
$G_\infty$-structure on $\g_2$ was proved to exist by 
Tamarkin~[Ta] and relies on a Drinfeld's associator. Tamarkin's $G_\infty$-morphism also restricts into a commutative, associative up to homotopy morphism ($C_\infty$-morphism for short). The $C_\infty$-structure on $\g_2$ (given by restriction of the $G_\infty$-one) highly depends on Drinfeld's associator, and any two choices of a Drinfeld associator yields {\it a priori} different $C_\infty$-structures. 
When $M$ is a Poisson manifold,
Kontsevich and Tamarkin homomorphisms imply the existence of a star-product (see
[BFFLS1]  and [BFFLS2] for a definition).
A connection between the two approaches has been given in
[KS] but the morphisms given by  Kontsevich and Tamarkin are not the same.
The aim of this paper is to show  that, given 
any $G_\infty$-structure on $\g_2$ and any  $C_\infty$-morphism $\varphi$ between
$\g_1$ and $\g_2$, there exists a $G_\infty$-morphism $\Phi$
between
$\g_1$ and $\g_2$ that restricts to $\varphi$. We also show that any $L_\infty$-morphism can be deformed into a $G_\infty$-one.

\smallskip

\noindent In the first section, we fix notation and recall the definitions of 
$L_\infty$ and $G_\infty$-structures. In the second section we state and prove 
the main theorem. In the last section we show that 
any two $G_\infty$-morphisms given by Tamarkin's method are homotopic.

\smallskip

\noindent {\bf Remark~:} In the sequel, unless otherwise is stated,  the manifold $M$ is supposed to be $\RM^n$ for some $n\geq 1$. Most results could be generalized to other manifolds using techniques of Kontsevich~[Ko1] (also see~[TS], [CFT]). 

\centerline {\bf 1-$C_\infty$, $L_\infty$ and $G_\infty$-structures}

\smallskip

\noindent For any graded vector space $\g$, 
we choose the following degree on $\ve^{\bullet}\g$~: if $X_1,\dots, X_k$ are homogeneous elements of respective degree $|X_1|,\dots |X_k|$, then
$$|X_1\wedge \cdots\wedge X_k|= |X_1|+\cdots+|X_k|-k.$$
In particular the component $\g=\ve^1\g \subset\ve^{\bullet}\g$ 
is the same as the space $\g$ with degree shifted by one. The space  $\ve^{\bullet}\g$ 
with the deconcatenation cobracket is the cofree cocommutative 
coalgebra on $\g$ with degree shifted by one (see [LS],
Section 2).
Any  degree 
one map $d^{k}:\ve^k\g \to \g$ ($k\geq 1$) extends into  
a derivation $d^k:\ve^{\bullet}\g\to \ve^{\bullet}\g$ 
of the coalgebra $\ve^{\bullet}\g$ by cofreeness property.

\proclaim{Definition 1.1}
A vector space $\g$ is endowed with a $L_\infty$-algebra 
(Lie algebras ``up to homotopy'') structure
if there are degree one linear  maps 
$m^{{1,\dots,1}}$, with $k$ ones~: $\ve^k \g \rightarrow \g$ such that if we 
extend them to  maps $\ve^\bullet \g \rightarrow \ve^\bullet \g$, then  $d \circ d=0$ where $d$ is the derivation 
$$d=m^1+m^{1,1}+\cdots+m^{1,\dots,1}+\cdots.$$
\endproclaim

\noindent For  more details  on $L_{\infty}$-structures, 
see~[LS]. 
It follows from the definition that a $L_\infty$-algebra 
structure induces a differential coalgebra structure 
on $\ve^\bullet \g$ and that the map  $m^1: \g \to \g$ 
is a differential. If $m^{1,\dots,1}$~: $\ve^k \g \rightarrow \g$
are $0$ for $k \geq 3$, we get the usual definition of
(differential if $m^{1}\not= 0$) graded Lie algebras.

\medskip

\noindent 
For any graded vector space $\g$, we 
denote $\underline{\g^{\ot n}}$ the quotient of $\g^{\ot n}$ 
by the image of all
shuffles of length $n$ (see [GK] or [GH] for
details). 
The graded vector space $\oplus_{n\geq 0}\underline{\g^{\ot n}}$ 
is a quotient coalgebra of the tensor coalgebra $\oplus_{n\geq 0}\g^{\ot n}$.  
It is well known  
that this coalgebra $\oplus_{n\geq 0}\underline{\g^{\ot n}}$ 
is the cofree Lie coalgebra on the vector space $\g$
(with degree  shifted by minus one).

\proclaim{Definition 1.2}
 A $C_\infty$-algebra 
(commutative and asssociative ``up to homotopy'' algebra) structure
on a vector space $\g$ is given by a collection of degree one linear  maps $m^{k}$~: $\underline{\g^{\otimes k}} \rightarrow \g$ such that if we 
extend them to  maps $\gh \rightarrow \gh$, then  $d \circ d=0$ where $d$ is the derivation 
$$d=m^1+m^{2}+m^{3}+\cdots.$$
\endproclaim
In particular a $C_\infty$-algebra is an $A_\infty$-algebra. 

\medskip

\noindent For any space $\g$, 
we denote $\ve^\bullet\underline{\g^{\otimes \bullet}}$ 
the graded space 
$$\ve^\bullet\underline{\g^{\otimes \bullet}}=\mathop{\oplus}\limits_{m \geq 1,~p_1+\cdots +p_n=m} 
\underline{\g^{\otimes p_1}}\wedge 
\cdots \wedge \underline{\g^{\otimes p_n}}.$$ 
We use the following grading on 
$\ve^\bullet\underline{\g^{\otimes \bullet}}$: 
for $x_1^1,\cdots, x_n^{p_n}\in\g$, we define
$$|\underline{x_1^1\otimes\cdots \otimes x_1^{p_1}}
\wedge \cdots\wedge \underline{x_n^1\otimes\cdots \otimes x_n^{p_n}}|=
\sum_{i_1}^{p_1}|x_1^{i_1}|+\cdots +\sum_{i_n}^{p_1}|x_n^{i_n}|-n. $$
Notice that the induced grading on
 $\ve^\bullet\g\subset \ve^\bullet\underline{\g^{\otimes \bullet}}$ 
 is the same than the one introduced above. 
The cobracket on $\oplus \underline{\g^{\otimes \bullet}}$  
and the coproduct on $\ve^\bullet \g$ extend to a cobracket 
and a coproduct on $\gq$ which yield a Gerstenhaber coalgebra structure on  
$\ve^\bullet\underline{\g^{\otimes \bullet}}$. 
It is well known that this coalgebra structure is cofree 
(see~[Gi], Section~3 for example).

\proclaim{Definition 1.3}
A $G_\infty$-algebra (Gerstenhaber algebra ``up to homotopy'') structure 
on a graded vector space $\g$ is given by a collection of degree one maps
$$m^{p_1,\dots,p_n}~: ~
\underline{\g^{\otimes p_1}}\wedge
\cdots \wedge \underline{\g^{\otimes p_n}} \rightarrow \g$$ indexed by $p_1,\dots p_n\geq 1$ such that 
their canonical  extension: 
$\ve^\bullet\underline{\g^{\otimes \bullet}}\rightarrow 
\ve^\bullet\underline{\g^{\otimes \bullet}}$
satisfies  $ d \circ d=0$ where $$ d=\sum\limits_{m \geq 1,~p_1+\cdots p_n=m}m^{p_1,\dots,p_n}.$$
\endproclaim

\noindent 
Again, as the coalgebra structure of $\ve^\bullet\underline{\g^{\otimes \bullet}}$ 
is cofree, the map $d$ makes $\ve^\bullet\underline{\g^{\otimes \bullet}}$  
into a differential coalgebra.
If the maps $m^{p_1,\dots,p_n}$ are $0$ for $(p_1,p_2,\dots)\not=(1,0,\dots)$,
$(1,1,0,\dots) $ or $(2,0,\dots)$, we get the usual definition of
(differential if $m^{1}\not= 0$) Gerstenhaber algebra.

\medskip

\noindent The space of multivector fields $\g_1$ is endowed with
a graded Lie  bracket $[-,-]_S$ called 
the Schouten bracket (see [Kos]).
This Lie algebra can be extended into a Gerstenhaber 
algebra, with commutative structure given by the exterior product:  
$(\alpha,\beta) \mapsto \alpha \ve \beta$

Setting $d_1=m_1^{1,1}+m_1^2$, where
$m_1^{1,1}$~: $\ve^2 \g_1 \rightarrow \g_1$, 
and 
$m_1^2$~: $\underline{\g_1^{\otimes 2}}\rightarrow \g_1$ are the extension
of the Schouten bracket and the exterior product, we find that $(\g_1,d_1)$ is a $G_\infty$-algebra.

\medskip

\noindent In the same way, one can define a differential Lie algebra structure 
on the vector 
space 
$\g_2=C(A,A)=\bigoplus_{k \geq 0}C^k(A,A)$,
the space of Hochschild cochains (generated by differential 
$k$-linear maps from
$A^k$ to $A$), where  $A=C^{\infty}(M)$ is the  algebra  of smooth differential
functions over $M$.
Its bracket $[-,-]_G$, called the Gerstenhaber bracket, is defined,
for $D,E \in \g_2$, by 
$$[D,E]_G=\{D | E\} - {(-1)}^{|E||D|}\{E| D\},$$where
$$\{D|E\}(x_1,{\dots},x_{d+e-1})=\sum_{i\geq
0}{(-1)}^{|E|{\cdot}i}D(x_1,{\dots},x_i,E(x_{i+1},{\dots},
x_{i+e}),{\dots}).$$
The space $\g_2$   has a grading defined by 
$\mid\! D\!\mid=k\Leftrightarrow D\in C^{k+1}(A,A)$ 
and its differential is $b=[m,-]_G$, 
where $m\in C^2(A,A)$ is the commutative multiplication on $A$. 

\smallskip

\noindent
Tamarkin (see [Ta] or also [GH]) 
stated the existence of a $G_\infty$-structure 
on $\g_2$ (depending on a choice of a Drinfeld associator) given by a differential 
$d_2=m_2^{1}+m_2^{1,1}+m_2^2+\cdots +
m_2^{p_1,\dots,p_n}+\cdots,$ on $\gd$ satisfying
$d_2 \circ d_2=0$. Although this structure is non-explicit, it satisfies the following three properties~:  
$$\eqalign{
(a)\hskip0.3cm& m_2^1 \hbox{ is the extension of the differential } b\cr
(b)\hskip0.3cm&m_2^{1,1}\hbox{ is the extension of the Gerstenhaber bracket }
[-,-]_G\cr
&\hbox{and } m_2^{1,1,\dots,1}=0 \cr
(c)\hskip0.3cm& m_2^2 \hbox{ induces the exterior product in cohomology and the} \cr
&\hbox{ collection of the }(m^k)_{k\geq 1}\hbox{ defines a }C_\infty\hbox{-structure on } \g_2.
\hfill  (1.1) 
}
$$

\medskip

\proclaim{Definition 1.4}
A $L_{\infty}$-morphism between two $L_\infty$-algebras 
$(\g_1,d_1=m^1_1+\dots)$ and $(\g_2,d_2=m^1_2+\dots)$ 
is  a morphism of differential coalgebras
$$
\varphi~: ~(\ve^\bullet \g_1,d_1) \rightarrow (\ve^\bullet \g_2,d_2).
\eqno (1.2) $$
\endproclaim
\noindent
Such a  map $\varphi$ is uniquely determined by a collection of  maps
$\varphi^n$~: $\ve^n \g_1 \rightarrow  \g_2$ (again by cofreeness properties).
In the case $\g_1$ and $\g_2$ are respectively the
graded Lie algebra $(\Gamma(M,\ve TM), ~[-,-]_S)$ and the differential
graded Lie algebra
$\left(C\left(A,A\right) ,~[-,-]_G\right)$,
the formality theorems of Kontsevich and Tamarkin
state the existence of a $L_\infty$-morphism between $\g_1$ and $\g_2$ such that $\varphi^1$ is the Hochschild-Kostant-Rosenberg 
quasi-isomorphism.

\medskip
\proclaim{Definition 1.5}
A morphism of $C_\infty$\!-algebras between two 
$C_\infty$\!ñ-algebras  $(\g_1,d_1)$ and $(\g_2,d_2)$ 
is a map $\phi : (\ghu,d_1) \to (\ghd,d_2)$ of codifferential coalgebras. 
\endproclaim
A $C_\infty$-morphism is in particular a morphism of $A_\infty$-algebras and is uniquely determined by maps $\p^k: \underline{\g^{\otimes k}}\to \g$.

\medskip

\proclaim{Definition 1.6}
A morphism of $G_\infty$\!-algebras between two 
$G_\infty$\!-algebras  $(\g_1,d_1)$ and $(\g_2,d_2)$ 
is a map $\phi : (\gu,d_1) \to (\gd,d_2)$ of codifferential coalgebras. 
\endproclaim
\noindent There are  coalgebras inclusions $\ve^{\bullet}\g \to \ve^\bullet\underline{\g^{\otimes \bullet}}$,  $\gh \to  \ve^\bullet\underline{\g^{\otimes \bullet}} $ and it is easy to check that any $G_\infty$-morphism between two $G_\infty$-algebras $\left(\g, \sum m^{p_1,\dots,p_n}\right)$, $\left({\g'}, 
\sum {m'}^{p_1,\dots,p_n}\right)$  restricts to a $L_\infty$-morphism $\left(\ve^{\bullet}\g, \sum m^{1,\dots,1}\right) \to \left(\ve^{\bullet}{\g'}, \sum {m'}^{1,\dots,1}\right)$ and a $C_\infty$-morphism $\left(\oplus\underline{\g^{\otimes \bullet}},\sum m^k\right)\to \left(\oplus\underline{{\g'}^{\otimes \bullet}},\sum {m'}^k\right)$.
In the case $\g_1$ and $\g_2$ are as above,
Tamarkin's theorem states that there exists a $G_\infty$-morphism between the two $G_\infty$ algebras $\g_1$ and $\g_2$ (with the
$G_\infty$ structure he built) that restricts to a $C_\infty$ and a $L_\infty$-morphism.

\medskip

\centerline {\bf 2-Main theorem}

\bigskip

\noindent We keep the notations of the previous section, in particular $\g_2$ is
the Hochschild complex of cochains on $C^\infty(M)$ and
$\g_1$ its cohomology. 
Here is our main theorem.

\proclaim{Theorem 2.1}
Given any $G_\infty$-structure $d_2$ on $\g_2$ satisfying the three properties of~(1.1),
and any $C_\infty$-morphism $\varphi$ between
$\g_1$ and $\g_2$ such that $\varphi^1$ is the Hochschild-Kostant-Rosenberg map, there exists a $G_\infty$-morphism $\Phi : (\g_1,d_1)\to (\g_2,d_2)$
 that restricts to $\varphi$. 

Also, given
 any $L_\infty$-morphism $\gamma$ between
$\g_1$ and $\g_2$ such that $\gamma^1$ is the Hochschild-Kostant-Rosenberg map, there exists a $G_\infty$-structure $(\g_1,d'_1)$ on $\g_1$ and $G_\infty$-morphism $\Gamma : (\g_1,d'_1)\to (\g_2,d_2)$
 that restricts to $\gamma$. Moreover there exists a $G_\infty$-morphism $\Gamma': (\g_1,d_1) \to (g_1,d_1)$.
\endproclaim
In particular, Theorem~2.1 applies to the formality map of Kontsevich and also to any $C_\infty$-map derived (see~[Ta], [GH]) from any $B_\infty$-structure on $\g_2$ lifting the Gerstenhaber structure of $\g_1$.

Let us first recall the proof of Tamarkin's formality theorem (see~[GH]
for more details): 
$$
\eqalign{
1.\hskip0.3cm& \hbox{First one proves  there exists
a }G_\infty\hbox{-structure on }\g_2,\hbox{ with differential }d_2,
\hbox{  as in~~(1.1)}. \cr
2.\hskip0.3cm& \hbox{Then, one constructs a }G_\infty\hbox{-structure on }\g_1\hbox{ given
by a differential }d_1'\hbox{ together with}\cr
&\hbox{a }ñG_\infty\hbox{-morphism
}\Phi\hbox{ between }(\g_1,d_1')\hbox{ and }(\g_2,d_2).\cr
3.\hskip0.3cm& \hbox{Finally, one constructs a }G_\infty\hbox{-morphism
}\Phi' \hbox{ between }(\g_1,d_1)\hbox{ and }(\g_1,d_1').
}
$$
The composition $\Phi \circ \Phi'$ is then a $G_\infty$-morphism between
$(\g_1,d_1)$ and $(\g_2,d_2)$,thus restricts to a $L_\infty$-morphism
between the differential graded Lie algebras $\g_1$ and $\g_2$.

\smallskip

\noindent We suppose  now that, in the first step,
we take any $G_\infty$-structure on $\g_2$ given by a differential $d_2$
and we suppose we are given a $C_\infty$-morphism $\varphi$ and a  $L_\infty$-morphism $\gamma$ between
 $\g_1$ and $\g_2$ satisfying $\gamma^1=\varphi^1=\varphi_{\HKR}$ the Hochschild-Kostant-Rosenberg quasi-isomorphism.

\smallskip

{ \noindent {\it Proof of Theorem~2.1: }

The Theorem  will follow if we prove that
steps 2 and 3 of Tamarkin's construction 
are still true with the extra conditions
that the restriction of the $G_\infty$-morphism $\Phi$ (resp. $\Phi'$)
on the $C_\infty$-structures is the $C_\infty$-morphism $\varphi:\g_1\to \g_2$ (resp. $\Id$).

\smallskip

Let us recall (see [GH]) that the constructions of $\Phi$ 
and $d_1'$ can be made by
induction.
For $i=1,2$ and $n\geq 0$, let us set 
$$V_i^{[n]}=\bigoplus_{p_1+\cdots +p_k=n}
~\underline{\g_i^{\otimes p_1}}\wedge
\cdots \wedge  \underline{\g_i^{\otimes p_k}}$$ 
and 
$V_i^{[\leq n]}=\sum_{k \leq n}V_i^{[k]}$.
Let $d_2^{[n]}$ and $d_2^{[\leq n]}$ be 
the sums 
$$d_2^{[n]}=\sum_{p_1+\cdots+p_k=n}d_2^{p_1,\dots,p_k} 
\qquad \hbox{and} \qquad d_2^{[\leq n]}=\sum_{p\leq n} d_2^{[p]}.$$
Clearly,  $d_2=\sum_{n\geq 1} d_2^{[n]}$.
In the same way, we denote $d_1'=\sum_{n \geq 1} {d'}_1^{[n]}$ with 
$${d'}_1^{[ n]}=\sum_{p_1+\cdots+p_k=n}d_1^{'p_1,\dots,p_k} \qquad \hbox{and} \qquad 
{d'}_1^{[\leq n]}=\sum_{1 \leq k \leq n} {d'}_1^{[k]}.$$
We know from Section~1
that a morphism $\Phi$~: $(\gu,d_1')\rightarrow
(\gd,d_2)$ is uniquely determined by its components 
$\Phi^{p_1,\dots,p_k}:\underline{\g_1^{\otimes p_1}}
\wedge \cdots \wedge \underline{\g_1^{\otimes p_k}}\rightarrow \g_2$.
Again, we have
$\Phi=\sum_{n \geq 1} \Phi^{[n]}$ with
$$\Phi^{[n]}=\sum_{p_1+\cdots+p_k=n}\Phi^{p_1,\dots,p_k}
\qquad \hbox{and} \qquad 
\Phi^{[\leq n]}=\sum_{1 \leq k \leq n} \Phi^{[k]}.$$
We want to construct the maps ${d'}_1^{[ n]}$ and $\Phi^{[n]}$ by induction
with the initial condition 
$${d'}_1^{[1]}=0 \qquad \hbox{and} \qquad \Phi^{[1]}=\varphi_{\HKR},$$
 where $\varphi_{\HKR}:(\g_1,0)\to (\g_2,b)$ is the 
Hochschild-Kostant-Rosenberg quasi-isomorphism 
(see~[HKR]) 
defined, for $\alpha \in \g_1$, $f_1,\cdots,f_n\in A$, by
$$\varphi_{\HKR}~:~\alpha \mapsto \big( (f_1,\dots, f_n)\mapsto
\; \langle\alpha, d f_1 \wedge \cdots \wedge d f_n\rangle \big).$$
Moreover, we want the following extra conditions to be true:
$$
\Phi^{k\geq 2}=\varphi^{k},\qquad 
{d'}_1^{2}=d_1^{2}, \qquad {d'}_1^{k\geq 3}=0.
\eqno (2.3)$$
Now suppose the construction is done for $n-1$ ($n\geq 2$), i.e., 
we have built maps $({d'}_1^{[i]})_{i \leq n-1}$ and $(\Phi^{[i]})_{i
\leq n-1}$ satisfying conditions (2.3) and
$$
\Phi^{[\leq n-1]} \circ {d'}_1^{[\leq n-1]}=d_2^{[\leq n-1]} \circ \
\Phi^{[\leq n-1]}\hbox{ on }V_1^{[\leq n-1]}
\hbox{ and }
{d'}_1^{[\leq n-1]} \circ {d'}_1^{[\leq n-1]}=0\hbox{ on }V_1^{[\leq n]}.
\eqno (2.4) $$
In [GH], we prove that for any such $({d'}_1^{[i]})_{i \leq n-1}$ 
and $(\Phi^{[i]})_{i
\leq n-1}$, one can construct ${d'}_1^{[n]}$ and $\Phi^{[n]}$
such that condition (2.4) is true for $n$ instead of $n-1$.
To complete the proof of Theorem~2.1 (step 2), we have to show that 
${d'}_1^{[n]}$ and $\Phi^{[n]}$ can be chosen to satisfy conditions
(2.3). In the equation~2.4, 
the terms
${d'}_1^k$ and $\Phi^k$ only act on $V_1^k$. 
So one can replace 
$\Phi^{n}$ with $\varphi^{n}$, 
${d'}_1^{2}$ with $d_1^{2}$ 
(or ${d'}_1^{i},i\geq 3$ with $0$)
provided conditions (2.4) are still satisfied on $V_1^{n}$. The
other terms acting on $V_1^{n}$ in the equation~(2.4) 
only involve terms $\Phi^{m}=\varphi^{m}$ and ${d'}_1^{m}$.
Then conditions~(2.4) on $V_1^{1,\dots, 1}$ 
are the equations that should be satisfied by a 
$C_\infty$-morphism between the $C_\infty$-algebras 
$(\g_1,{d'}_1^{1,1}=d_1^{1,1})$ and $(\g_2,\sum_{k\geq 1} d_2^{k})$ 
restricted to $V_1^{n}$. Hence by hypothesis on $\varphi$ 
the conditions hold. 

\smallskip

Similarly the construction of $\Phi'$ 
can be made by induction. Let us recall the proof given in
[GH].
Again a morphism $\Phi'$~: $(\gu,d_1)\rightarrow
(\gd,d_1')$ is uniquely determined by its components 
${\Phi'}^{p_1,\dots,p_k}:\underline{\g_1^{\otimes p_1}}
\wedge \cdots \wedge \underline{\g_1^{\otimes p_k}}\rightarrow \g_1$.
We write
$\Phi'=\sum_{n \geq 1} {\Phi'}^{[n]}$ with
$${\Phi'}^{[n]}=\sum_{p_1+\cdots+p_k=n}{\Phi'}^{p_1,\dots,p_k}
\qquad \hbox{and} \qquad 
{\Phi'}^{[\leq n]}=\sum_{1 \leq k \leq n}{\Phi'}^{[k]}.$$
We construct the maps ${\Phi'}^{[n]}$ by induction
with the initial condition ${\Phi'}^{[1]}=\Id$.
Moreover, we want the following extra conditions to be true:
$$
{\Phi'}^{n}=0 \hbox{ for } n \geq 2. 
\eqno (2.5)$$
Now suppose the construction is done for $n-1$ ($n\geq 2$), i.e., 
we have built maps $({\Phi'}^{[i]})_{i
\leq n-1}$ satisfying conditions (2.5) and
$$
{\Phi'}^{[\leq n-1]} {d}_1^{[\leq n]}={d_1'}^{[\leq n]} 
{\Phi'}^{[\leq n-1]} \hbox{ on } V_1^{[\leq n]}.
\eqno (2.6) $$
In [GH], we prove that for any such $({\Phi'}^{[i]})_{i
\leq n-1}$, one can construct ${\Phi'}^{[n]}$
such that condition (2.6) is true for $n$ instead of $n-1$
in the following way~: 
making the equation ${\Phi'}  {d}_1={d_1'}  \Phi'$
on $V_1^{[n+1]}$ explicit, we get
$$
{\Phi'}^{[\leq n]}\, {d}_1^{[\leq n+1]}={d_1'}^{[\leq n+1]}\,{\Phi'}^{[\leq
n]}.
\eqno (2.7)$$
If we now take into account that $d_1^{[i]}=0$ for $i\not=2$,
${d_1'}^{[1]}=0$
and that on $V_1^{[n+1]}$ we have ${\Phi'}^{[ k]} {d}_1^{[l]}=
{d_1'}^{[\leq k]}{\Phi'}^{[l]}=0$ for $k+l > n+2$, the identity (2.7)
becomes
$${\Phi'}^{[\leq n]}\,d_1^{[2]}=\sum_{k=2}^{n+1} {d_1'}^{[k]}\,{\Phi'}^{[\leq
n-k+2]}.$$
As
${d_1'}^{[2]}={d_1}^{[2]}$, (2.7) is equivalent to
$${d_1}^{[2]}{\Psi'}^{[\leq n]}-{\Phi'}^{[\leq n]}{d_1}^{[2]}=
\left[{d_1}^{[2]},{\Phi'}^{[\leq n]}\right]=
-\sum_{k=3}^{n+1} {d_1'}^{[k]}{\Phi'}^{[\leq
n-k+2]}.$$
Notice that  $d_1^{[2]}=m_1^{1,1}+m_1^2$. Then the  construction will be possible 
when the term $\sum_{k=3}^{n+1} {d_1'}^{[k]}{\psi'}^{[\leq
n-k+2]}$ is a couboundary in the subcomplex of $({\hbox{End}}(\gu),$ $[d_1^{[2]},-])$ consisting of maps which restrict to zero on $\oplus_{n\geq2}\underline{\g_1}^{\otimes n}$. It is always a cocycle by straightforward computation (see~[GH]) and the subcomplex is acyclic because both  $({\hbox{End}}(\gu),$ $[d_1^{[2]},-])$ and the Harrison cohomology of $\g_1$ are trivial according to Tamarkin~[Ta] (see also [GH] Proposition 5.1 and [Hi] 5.4). 

In the case of the $L_\infty$-morphism $\gamma$, the first step is similar: the fact that $\gamma$ is a $L_\infty$-map enables us to build a $G_\infty$-structure  $(\g_1,d'_1)$ on $\g_1$ and a $G_\infty$-morphism $\Gamma:(\g_1,d'_1)\to (\g_2,d_2)$ such that:
$$
\Gamma^{1,\dots,1}=\gamma^{1,\dots,1},\qquad 
{d'}_1^{1,1}=d_1^{1,1},\qquad {d'}_1^{1,1,\dots,1}=0.
\eqno (2.8)$$
For the second step, we have to build a map $\Gamma'$  satisfying the equation
$${d_1}^{[2]}{\Gamma'}^{[\leq n]}-{\Gamma'}^{[\leq n]}{d_1}^{[2]}=
\left[{d_1}^{[2]},{\Gamma'}^{[\leq n]}\right]=
-\sum_{k=3}^{n+1} {d_1'}^{[k]}{\Gamma'}^{[\leq
n-k+2]}$$ on $V_1^{[n+1]}$ for any $n\geq 1$. Again,  because Tamarkin has prooved that the complex $({\hbox{End}}(\gu),$ $[d_1^{[2]},-])$ is acyclic (we are in the case
$M=\RM^n$), the result follows from the fact that $\sum_{k=3}^{n+1} {d_1'}^{[k]}{\Gamma'}^{[\leq
n-k+2]}$ is a cocycle. The difference with the $C_\infty$-case is that the  ${\Gamma'}^{1,\dots,1}$ could be non zero. 
{\null \hfill  $\blacksquare$ \par\medskip}

\medskip

\centerline {\bf 3-The difference between two $G_\infty$-maps}

\bigskip

In this section we investigate the difference between 
 two differents $G_\infty$-formality maps.

\smallskip

\noindent We fix once for all a $G_\infty$-structure 
on $\g_2$ (given by a differential $d_2$) 
satisfying the conditions~(1.1) and a morphism 
of $G_\infty$-algebras $T:(\g_1, d_1) \to (\g_2,d_2)$ 
such that $T^1: \g_1\to \g_2$ is $\varphi_{\HKR}$. 
Let $K:(\g_1, d_1) \to (\g_2,d_2)$ be any other  
$G_\infty$-morphism with $K^1=\varphi_{\HKR}$ 
(for example any lift of a Kontsevich formality map or any $G_\infty$-maps lifting another $C_\infty$-morphism).

\proclaim{Theorem 3.1}
There exists a map $h:\gu\to \gd$ such that $$ T-K=h\circ d_1+d_2\circ h.$$
\endproclaim

In other words the formality $G_\infty$-morphisms $K$ and $T$ are homotopic. 

\smallskip

\noindent The maps $T$ and $K$ are elements of the cochain complex $\left(\Hom(\gu,\gd), \delta\right)$  with  differential given, for all $f\in \Hom(\gu,\gd), |f|=k$, by 
$$\delta(f)= d_2\circ f-(-1)^{k}f\circ d_1 .$$

\smallskip

\noindent We first compare this cochain complex with 
the complexes 
$\left(\End(\gu),[d_1;-]\right)$ and 
$\left(\End(\gd),[d_2;-]\right)$ (where $[-;-]$ is the graded 
commutator of morphisms). There are morphisms 
$$T_*:\End(\gu)\to \Hom(\gu,\gd), \ T^*:  \End(\gd)\to \Hom(\gu,\gd)$$ defined, for $f\in  \End(\gd)$ and $g\in \Hom(\gu,\gd)$, by 
$$T_*(f)=T\circ f,\ T^*(g)= g\circ T.$$ 

\proclaim{Lemma 3.2}The morphisms 
$$T_*:\left(\End(\gu),[d_1;-]\right)\to \left(\Hom(\gu,\gd),\delta\right)\leftarrow \left( \End(\gd),[d_2;-]\right):T^*$$ of cochain complexes are quasi-isomorphisms.
\endproclaim
\noindent Remark: This lemma  holds for every manifold $M$ and any $G_\infty$-morphism $T:(\g_1,d_1)\to (g_2,d_2)$. 

\smallskip

\demo{Proof :}
First we show that $T_*$ is a morphism of complexes. Let $f\in \End(\gd)$ with $|f|=k$, then
$$
\eqalign{
T_*([d_1;f])=& T\circ d_1\circ f -(-1)^kT\circ f\circ d_1\cr
=& d_2\circ(T\circ f)-(-1)^k(T\circ f)\circ d_1\cr
=& \delta (T_*(f)).
}
$$
Let us prove now that $T_*$ is a quasi-isomorphism. 
For any graded vector space $\g$, the space  $\gq$ has the structure of a  
filtered space where the $m$-level of the filtration is 
$F^m(\gq)=\oplus_{p_1+\dots+p_n-1\leq m}
\underline{\g^{\otimes p_1}}\wedge \dots \underline{\g^{\otimes p_n}}$. 
Clearly the differential $d_1$ and $d_2$ are compatible with the filtrations 
on $\gu$ and $\gd$, hence $\End(\gu\!)$ and $\Hom(\gu\! ,\gd\!)$ 
are filtered cochain complex. This yields two spectral sequences 
(lying in the first quadrant) $E_\bullet^{\bullet,\bullet}$ and 
$\widetilde{E}_{\bullet}^{\bullet, \bullet}$ which converge respectively 
toward the cohomology $H^{\bullet}(\End(\gu))$ and 
$H^\bullet (\Hom(\gu,\gd))$. By standard spectral sequence techniques
it is enough to prove that the map 
$T^0_*~:~E_0^{\bullet, \bullet}\to \widetilde{E}_0^{\bullet, \bullet}$ induced by 
$T_*$ on the associated graded is a quasi-isomorphism.

\smallskip

\noindent  The induced differentials on $E_0^{\bullet, \bullet}$ 
and $\widetilde{E}_0^{\bullet, \bullet}$ are respectively 
$[d_1^1,-]=0$ and $d_2^1\circ (-)-(-)\circ d_1^1=b\circ(-)$ where 
$b$ is the Hochschild coboundary. 
By cofreeness property
we have the following two isomorphisms
$$E_0^{\bullet, \bullet}\cong  \End(\g_1), \qquad \widetilde{E}_0^{\bullet, \bullet} \cong \Hom(\g_1,\g_2).$$
The map $T_*^0:E_0^{\bullet, \bullet}\to \widetilde{E}^0_{\bullet \bullet}$ 
induced by $T_*$  is $\varphi_{\HKR}\circ(-)$. 
Let $p:\g_2\to \g_1$ be the projection onto the cohomology, 
 {\it i.e.} $p\circ \varphi_{\HKR}= \Id$.  
Let $u:\g_1\to \g_2$ be any map satisfying $b(u)=0$ and set 
$v=p\circ u \in \End(\g_1)$. 
One can choose a map $w:\g_1\to \g_2$ which satisfies for any 
$x\in \g_1$ the following identity  
$$\varphi_{\HKR}\circ p \circ u(x)-  u(x)=b\circ w(x).$$ 
It follows that $\varphi_{\HKR}(v)$ has the same class of 
homology as $u$ which proves the surjectivity of  $T_*^0$ in cohomology.  
The identity $p\circ \varphi_{\HKR}= \Id$ implies easily that 
$T_*^0$ is also injective in cohomology which finish the proof of the 
lemma for $T_*$.

\smallskip

\noindent The proof that $T^*$ is also a quasi-isomorphism is  analogous.
\enddemo

\smallskip

\demo{Proof of Theorem~3.1:}

\noindent It is easy to check that $T-K$ is a cocycle in 
$\left(\Hom(\gu,\gd),\delta \right)$.
The complex of cochain $\left(\End(\gu),[d_1,-]\right)\cong 
\left(\Hom(\gu,\g_1),[d_1,-]\right)$ is trigraded with $|\,|_1$ 
being the degree coming from the graduation of $\g_1$ and any 
element $x$ lying in $\underline{\g_1^{\otimes p_1}}\wedge \dots 
\wedge \underline{\g_1^{\otimes p_q}}$ satisfies $|x|_2=q-1$, 
$|x|_3=p_1+\dots p_q-q$.  
In the case $M=\RM^n$,  the cohomology $H^\bullet\left(\End(\gu),[d_1,-]\right)$ 
is concentrated in bidegree $(|\,|_2,|\,|_3)=(0,0)$  
(see~[Ta], [Hi]).  
By Lemma~3.2}, this is also the case for the cochain 
complex $\left(\Hom(\gu,\gd),\delta \right)$. 
Thus, its cohomology classes are determined by complex morphisms $(\g_1,0)\to
(\g_2, d_2^1)$ and  it is enough to prove that $T$ and $K$ determine the same
complex morphism $(\g_1,0)\to (\g_2, d_2^1=b)$ which is clear because $T^1$ and $K^1$ are both equal to the Hochschild-Kostant-Rosenberg map. 
{\null \hfill  $\blacksquare$ \par\medskip}

\smallskip

\proclaim{Remark}
It is possible to have an explicit formula for the map $h$ 
in Theorem~.3.1. In fact the quasi-isomorphism coming 
from Lemma~3.2 can  be made explicit using  explicit homotopy 
formulae for the Hochschild-Kostant-Rosenberg map (see~[Ha] 
for example) and deformation retract techniques 
(instead of spectral sequences) as in~[Ka]. 
The same techniques also apply to give explicit formulae for the
quasi-isomorphism giving the acyclicity of $\left(\End(\gu),[d_1;-]\right)$ 
in the proof of
theorem~3.1 (see~[GH] for example)
\endproclaim

\vskip15pt

\enddocument